\documentclass[12pt,epsf]{article}
\usepackage{amssymb,amsmath,eucal,amsthm}
\usepackage{graphicx}

\def\ds{\displaystyle}
\def\pa{\partial}
\def\bce\{\begin{center}
\def\ece{\end{center}}

\theoremstyle{definition}

\theoremstyle{remark}

\renewcommand{\proofname}{\bf Proof}
\newcommand{\N}{{I \!\! N}}
\newcommand{\R}{{I \!\! R}}
\newcommand{\C}{{I \!\! C}}
\newcommand{\ol}{\overline}

\title{\bf Stochastic generalized fractional HP equations and applications.}

\author{\bf I. D. Albu$^{a}$, M. Neam\c tu$^{b}$\thanks{Corresponding author}, D. Opri\c s$^{c}$}
\date{ }

\begin{document}
\maketitle

\begin{tabular}{cccccccc}
\scriptsize{$^{a}$ {\it Department of
Mathematics, Faculty of Mathematics,}}\\
\scriptsize{{\it West University of Timi\c soara, Bd. V. Parvan, nr. 4, 300223, Timi\c soara, Romania,}}\\
\scriptsize{{\it E-mail: albud@math.uvt.ro}}\\
\scriptsize{$^{b}${\it Department of Economic Informatics, Mathematics and Statistics, Faculty of Economics,}}\\
\scriptsize{{\it West University of Timi\c soara, Pestalozzi Street, nr. 16A, 300115, Timi\c soara, Romania,}}\\
\scriptsize{{\it E-mail:mihaela.neamtu@feaa.uvt.ro,}}\\
\scriptsize{$^{c}$ {\it Department of Applied
Mathematics, Faculty of Mathematics,}}\\
\scriptsize{{\it West University of Timi\c soara, Bd. V. Parvan, nr. 4, 300223, Timi\c soara, Romania,}}\\
\scriptsize{{\it E-mail: opris@math.uvt.ro}}\\

\end{tabular}
\begin{abstract} In this paper we established the condition for a curve to satisfy stochastic
generalized  fractional HP (Hamilton-Pontryagin) equations. These
equations are described using It\^{o} integral. We have also
considered the case of stochastic generalized  fractional
Hamiltonian equations,
 for a hyperregular Lagrange function. From the stochastic ge\-ne\-ra\-lized fractional Hamiltonian
  equations, Langevin generalized fractional equations were found and numerical simulations were done.

\vspace{0.2cm}

\noindent{\small {\bf AMS subject classification}: 34K50, 35L65,
26A33, 37N99, 60H10, 65C20, 65C30.}

\noindent{\small {\bf Keywords}: HP equations, stochastic
generalized  fractional Hamilton equations, hyperregular Lagrange
function, generalized fractional Langevin equations, Euler scheme.}

\end{abstract}

\section{Introduction}

\hspace{0.5cm} J.M. Bismut was the first one who introduced concepts
of stochastic geometric mechanics, in his article from 1981, when he
defined the notion of "stochastic Hamiltonian system". He showed
that that the stochastic flow of certain randomly perturbed
Hamiltonian systems on flat spaces ex\-tre\-mizes a stochastic
action, and using this property, he proved symplecticity and the
Noether theorem for stochastic Hamiltonian systems. Since then,
there has been a need in to find out tools and algorithms for the
study of this kind of systems with uncertainty. Bismut's work was
continued by Lazaro-Cami and Ortega ([11], [12]), in the sense that
his work was generalized to ma\-ni\-folds. Stochastic Hamiltonian
systems on manifolds extremize a stochastic action on the space of
manifold valued semimartingales, the reduction of stochastic
Hamiltonian system on the cotangent bundle of a Lie group, a counter
example for the converse of Bismut's original theorem.
    Very important in many scientific domains is fractional calculus: fractional derivatives,
fractional integrals, of any real or complex order. Fractional
calculus is used when fractional integration is needed. It is used
for studying simple dynamical systems, but it also describes complex
physical systems. For example, applications of the fractional
calculus can be found in chaotic dynamics, control theory,
stochastic modelling, but also in finance, hydrology, biophysics,
physics, astrophysics, cosmology, economics and so on ([2], [4],
[5], [9], [10]). But some other fields have just started to study
problems from fractional point of view. It is very fashionable to
study the fractional problems of the calculus of variations and
Euler-Lagrange type equations. The most famous fractional integrals
are Riemann-Liouville, Caputo, Grunwald-Letnikov and the most
frequently used is the Riemann-Liouville fractional integral. The
study of Euler-Lagrange fractional equations was continued by
Agrawal ([1], [6], [8]) that described these equations using the
left, respectively right fractional derivatives in the
Riemann-Liouville sense. Standard multi-variable variational
calculus also has some limitations. But in [13], C.Udriste and D.
Opris showed that these limitations can by broken using the
multi-linear control theory. In [7] the novel concepts of fractional
action-like variational approach (FALVA) with time-dependent
fractional exponent and exponential time-dependent term is
introduced.
 In this paper, we restrict our attention to stochastic generalized
 fractional Hamiltonian systems characterized by Wiener processes
  and assume that the space of admissible curves in
configuration space is of class $C^1$. Random effects appear in the
balance of momentum equations, as white noise, that is why we may
consider randomly perturbed mechanical systems. It should be
mentioned that the ideas in this paper can be readily extended to
stochastic Hamiltonian systems driven by more general
semimartingales, but for the sake of clarity we restrict ourselves
to Wiener processes.
  In this paper we use the generalized left fractional Riemann-Liouville
   integral defined as a mixture of the fractal  action from physics and
    the discounted action at rate $\rho $, given in [7].
 Within this context, the results of the  paper are as follows:

1. The paper presents the results from [3] which show that almost
surely that a curve satisfies stochastic  HP equations if and only
if it extremizes a stochastic action. Suggestive examples and
numerical simulations are done.

2. Generalized fractional HP equations are described using the
generalized fractional Riemann-Liouville integral and the fractional
It\^{o} integral;

3. Langevian type stochastic generalized fractional equations are
obtained in the case of a hyperregular Lagrange function. Relevant
examples and numerical simulations are presented.

The paper is organized as follows: In Section 2, Hamilton-Pontryagin
(HP) principle is given to the stochastic setting to prove that a
class of mechanical systems with multiplicative noise appearing as
forces and torques possess a variational structure. For a
hyperregular Lagrange function, we get the stochastic  Hamiltonian
equations that lead to Langevin equations. Examples and numerical
simulations for the Lagrangian which describes the Samuelson model
from economics [5] are given. In Section 3, we extend the
generalized fractional Hamilton-Pontryagin (HP) principle to the
stochastic setting to prove that a class of mechanical systems with
multiplicative noise appearing as forces and torques possesses a
variational structure. For a hyperregular Lagrange function, we
obtain the stochastic generalized fractional Hamiltonian equations
that lead to Langevin generalized fractional equations. For a
Lagrange function, defined on $\R^2$, the corresponding generalized
fractional Langevin equations are simulated.  The generalized
fractional Hamiltonian and the Lagrangian description are joined
together to get the generalized fractional HP system.

\section{Stochastic HP mechanics}

\hspace{0.5cm} In this section a variational principe is introduced
for a class of stochastic Hamilton systems on manifolds. The
stochastic action is a sum of the classical action and stochastic
integral. The key feature of this principle is that one can recover
stochastic Hamilton equations for these systems. Roughly speaking,
this is accomplished by means of taking variations of this action
within the space of curve only (not the probability space) and
imposing the condition that this partial differential of the action
must be zero.

Let a paracompact, configuration manifold $Q$ and
$J^1(\R,Q)=\R\times TQ$, $T^\ast Q$, the associate bundle of $Q$.
Let $(\Omega, {\cal {P}}, P)$ be a probability space and $(w(t),
{\cal F}_t)_{t\in [a,b]}$, where $[a,b]\subset \R$, $w(t)$ is a
real-valued Wiener process and ${\cal F}_t$ is the filtration
generated by the Wiener process [3].

The paper adopts an HP viewpoint to develop a Lagrangian description
of stochastic Hamiltonian systems [3]. The HP principle unifies the
Hamiltonian and Lagrangian descriptions of mechanical system. The
classical HP action integral will be perturbed using deterministic
function $\gamma :Q\to\R$.

We consider the Lagrangian ${\cal L}:J^1(\R, Q)\to\R $. In the
stochastic context the HP principle states the following critical
point condition on $J^1(\R, Q)\oplus T^\ast Q$ for stochastic HP
action integral given by [3]:

\begin{equation}
\begin{split}
{\cal A}(q,v,p)\!=\!\int_a^b[{\cal L}
(s,q(s),v(s))\!+\!<p(s),\ds\frac{dq}{ds}\!-\!v(s)>]ds\!+\!
\int_a^b\gamma (q(s))dw(s)
\end{split}
\end{equation} where
\begin{equation*}
\begin{split}
& {\cal A}:\Omega\times C(PQ)\to\R\\
& C(PQ)=\{(s,q,v,p)\in C^0([a,b], PQ)|q\in C^1([a,b], Q), q(a)=q_a,
q(b)=q_b\},
\end{split}
\end{equation*} $[a,b]\subset \R$, $q_a, q_b\in Q$.

The action integral in the above principle consists of one Lebeque
integral with respect to s and It$\hat{o}$ stochastic integral with
respect to w. The action is random; i.e. for every sample point
$\omega\in\Omega$ we will obtain a different time-dependent
Lagrangian system. We will use the following notation $q(\omega,
s)=q(s)$, $v(\omega, s)=v(s)$, $p(\omega, s)=p(s)$. The HP path
space is a smooth infinite dimensional manifold. One can show that
is tangent space at $c=(q,v,p)\in C([a,b], q_1,q_2)$ consists of
maps $w=(q,v,p,\delta q, \delta v, \delta p)$ $\in C^0([a,b],
T(PQ))$ such that $\delta q(a)=\delta q(b)=0$ and $q, \delta q$ are
of class $C^1$. Let $(q,v,p)(\cdot , \varepsilon)\in C(PQ)$ denote a
one-parameter family of curves in ${\cal C}$ that is differential
with respect to $\varepsilon$. Define the differential of ${\cal A}$
as

\begin{equation*}
d{\cal A}(\delta q, \delta v, \delta p)=\ds\frac{\partial }{\partial
\varepsilon}{\cal A}(\omega , q(s, \varepsilon), v(s, \varepsilon),
p(s, \varepsilon))|_{\varepsilon=0}
\end{equation*}
where
\begin{equation*}
\begin{split}
& \delta q(s)=\ds\frac{\partial }{\partial \varepsilon}q(s,
\varepsilon)|_{\varepsilon=0}, \delta q(a)=\delta q(b)=0,\\
& \delta v(s)=\ds\frac{\partial }{\partial \varepsilon}v(s,
\varepsilon)|_{\varepsilon=0}, \delta p(s)=\ds\frac{\partial
}{\partial \varepsilon}p(s, \varepsilon)|_{\varepsilon =0}.
\end{split}
\end{equation*}

In terms of this differential one can state the following critical
point condition:

{\bf Theorem 1. [3]} {\it Let ${\cal L}:J^1(\R,Q)\to\R$ be a
Lagrangian on $J^1(\R,Q)$ of class $C^2$ with respect to $t, q, v$
and with the globally Lipschitz first derivative with respect to $t,
q$ and $v$. Let $\gamma :Q\to\R$ be a class $C^2$ function and with
the globally Lipschitz first derivative. Then almost certainly a
curve $c=(q,v,p)\in C(PQ)$ satisfies the stochastic HP equations:
\begin{equation}
\begin{split}
& dq^i=v^ids,\\
& dp_i=\ds\frac{\partial {\cal L}}{\partial q^i}ds+\ds\frac{\partial
\gamma}{\partial q^i}dw(s)\\
& p_i=\ds\frac{\partial {\cal L}}{\partial v^i}, i=1..n
\end{split}
\end{equation} if and only if it is a critical point of the function
${\cal A}:\Omega\times {\cal C}(PQ)\to\R$, i.e. d{\cal A}(c)=0.}

Let ${\cal L}:J^1(\R, Q)\to\R$ be a Lagrangian on $J^1(\R, Q)$,
hyperregular, that means $det\left ( \ds\frac{\pa^2{\cal L}}{\pa
v^i\pa v^j}\right )\neq 0$.

From (2) the following propositions are obtained:

{\bf Proposition 1. (Stochastic Hamilton equations).}  {\it The
equations (2) are equivalent to the following equations:
\begin{equation}
\begin{split}
& dq^i=\ds\frac{\pa H}{\pa p_i}ds,\\
& dp_i=-\ds\frac{\partial H}{\partial q^i}ds+\ds\frac{\partial
\gamma}{\partial q^i}dw(s), i=1..n
\end{split}
\end{equation} where $H=p_iv^i-L(t,q,v)$.}

The equations (3) represent Lagevin equations.

{\bf Proposition 2. } {\it If ${\cal L}=\ds\frac{1}{2}g_{ij}v^iv^j$,
where $g_{ij}$ are the components of a metric on the manifold $Q$,
equations (2) take the form:
\begin{equation}
\begin{split}
& dq^i=v^idt, dp_i=-g_{ij}v^jds+\ds\frac{\partial \gamma
(q)}{\partial
q^i}dw(s),\\
& dv^i=-\Gamma ^i_{jk}v^jv^kds+g^{ij}\ds\frac{\pa v(q)}{\pa
q^j}dw(s), i,j=1..n,,
\end{split}
\end{equation} where $\Gamma^i_{jk}$ are Cristoffel coefficients
associated to the considered metric.

The equations (3) become:
\begin{equation}
\begin{split}
& dq^i=g^{ij}p_jds, \\
& dp_i=\ds\frac{1}{2}\ds\frac{\pa g_{kl}}{\pa
q^i}p^kp^lds+\ds\frac{\partial \gamma (q)}{\partial q^i}dw(s).
\end{split}
\end{equation} }

{\bf Proposition 3.} {\it If ${\cal L}:J^1(\R, \R^n)\to\R$ is given
by:
\begin{equation}
{\cal L}=\ds\frac{1}{2}\delta_{ij}v^iv^j-V(q), q\in\R^n,
\end{equation} the equations (2) take the form:
\begin{equation}
\begin{split}
& dq^i=v^{i}ds, \\
& dp_i=-\ds\frac{\pa V}{\pa q^i}ds+\ds\frac{\partial \gamma
}{\partial q^i}dw(s), p_i=\delta_{ij}v^j.
\end{split}
\end{equation} The equations (3) become}:
\begin{equation}
\begin{split}
& dq^i=\delta^{ij}p_jds, \\
& dp_i=-\ds\frac{\pa V(q)}{\pa q^i}ds+\ds\frac{\partial \gamma
(q)}{\partial q^i}dw(s).
\end{split}
\end{equation}

{\bf Proposition 4.} {\it If ${\cal L}:J^1(\R, \R^n)\to\R$ is given
by:
\begin{equation}
{\cal L}=e^{-\rho s}L(q,v), q\in\R^n
\end{equation} the equations (2) take the form:
\begin{equation}
\begin{split}
& dq^i=v^{i}ds, \\
& dp_i=e^{-\rho s}\ds\frac{\pa L}{\pa q^i}ds+\ds\frac{\partial
\gamma (q)}{\partial q^i}dw(s),\\
& p_i=e^{-\rho s}\ds\frac{\pa L}{\pa v^i}.
\end{split}
\end{equation}

{\bf Proposition 5. (Samuelson. [5])} {\it If ${\cal L}:J^1(\R,
\R^2)\to\R$ is given by:
\begin{equation}
{\cal L}=-\ds\frac{1}{2}e^{-\rho s}(v^2+2avq+q^2), q\in\R,
\end{equation} the equations (2) take the form}:
\begin{equation}
\begin{split}
& dq=vds, \\
& dp=-e^{-\rho s}(av+q)ds+\ds\frac{\partial \gamma (q)}{\partial
q^i}dw(s),\\ p=-e^{-\rho s}(v+aq).
\end{split}
\end{equation} The equations (12) become}:
\begin{equation}
\begin{split}
& dq=-(aq+e^{\rho s}p)ds, \\
& dp=((a^2-1)e^{-\rho s}q+ap)ds+\ds\frac{d\gamma }{dq}.
\end{split}
\end{equation}

If $V(q)=q^2$, the Euler scheme for (13) is:
\begin{equation}
\begin{split}
& q(n+1)=q(n)-h(aq(n)+e^{\rho n}p(n)) \\
& p(n+1)=p(n)+h((a^2-1)e^{-\rho n}q(n)+ap(n))+q(n)G(n), n=0..N-1,
\end{split}
\end{equation} where $T>0$, $h=\ds\frac{T}{N}$, $N>0$,
$G(n)=w((n+1)h)-w(nh)$ and $q(n)=q(\omega, nh)$, $p(n)=p(\omega,
nh)$, $\rho\in (0,1)$, $a\in (-1,1)$.

For $\rho =0.003$, $a=0.03$, $h=0.001$ using Maple 13, the orbit
$(n, q(nh))$ is represented in Fig 1 and $(n, q(\omega, nh))$ in Fig
2:

\begin{center}\begin{tabular}{cc}
\\ Fig 1. $(n, q(nh))$ & Fig 2. $(n, q(\omega, nh))$\\
\includegraphics[width=5cm]{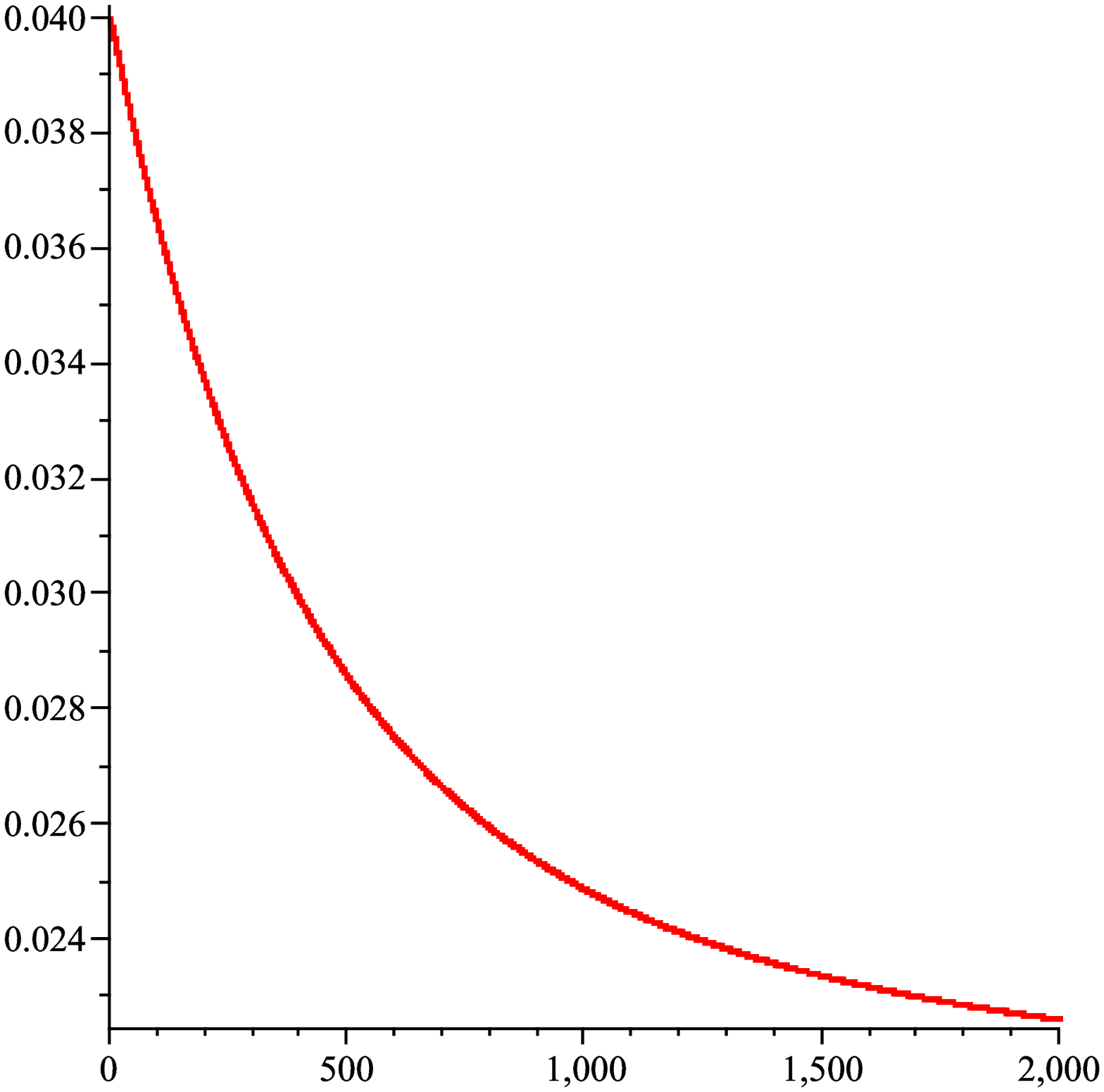} & \includegraphics[width=5cm]{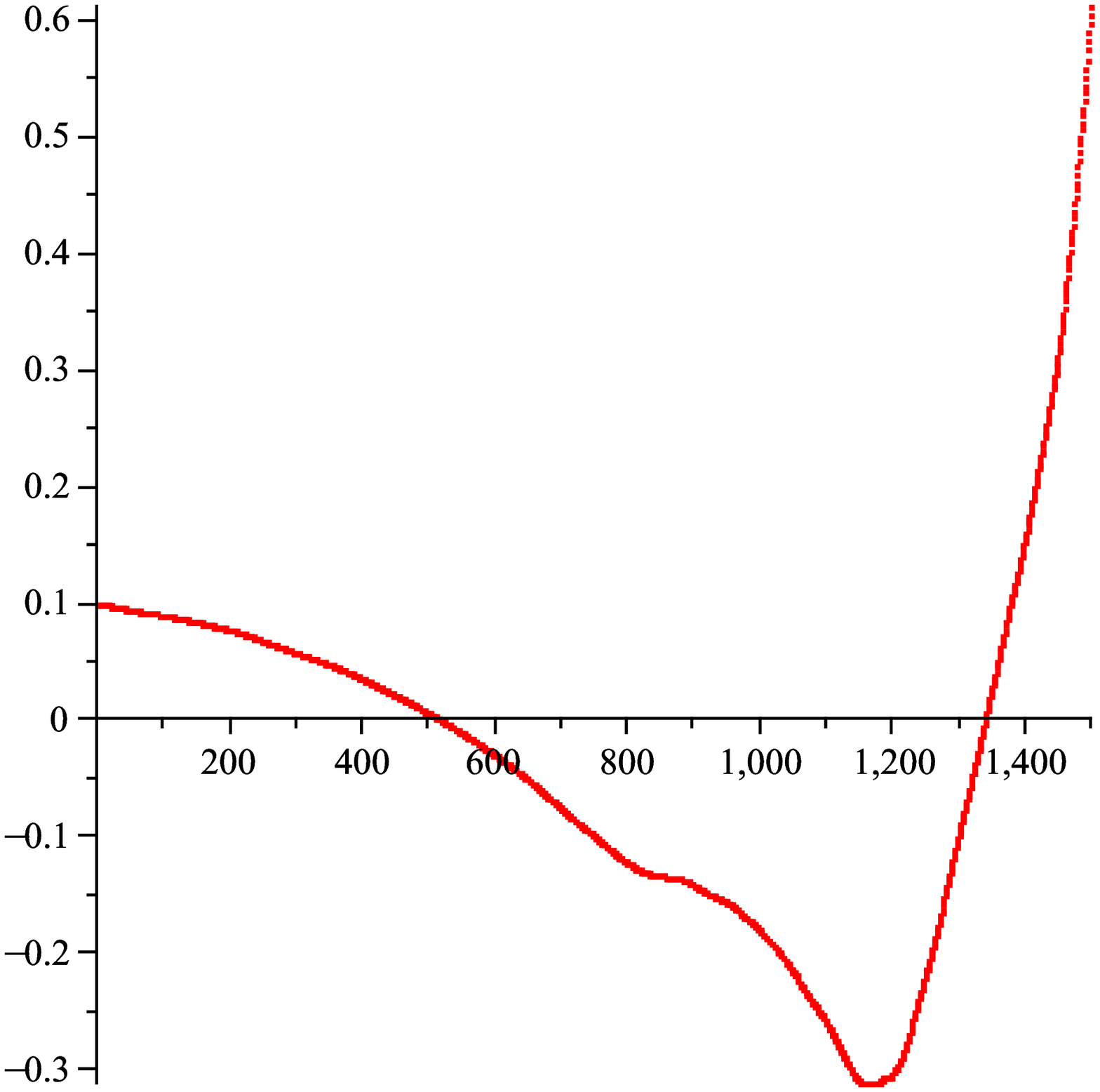}
\end{tabular}
\end{center}

In Figures 3 and 4 we can visualize the orbits $(n, p(nh))$, $(n,
p(\omega, nh))$:

\begin{center}\begin{tabular}{cc}
\\ Fig 3. $(n, p(nh))$ & Fig 4. $(n, p(\omega, nh))$\\
\includegraphics[width=5cm]{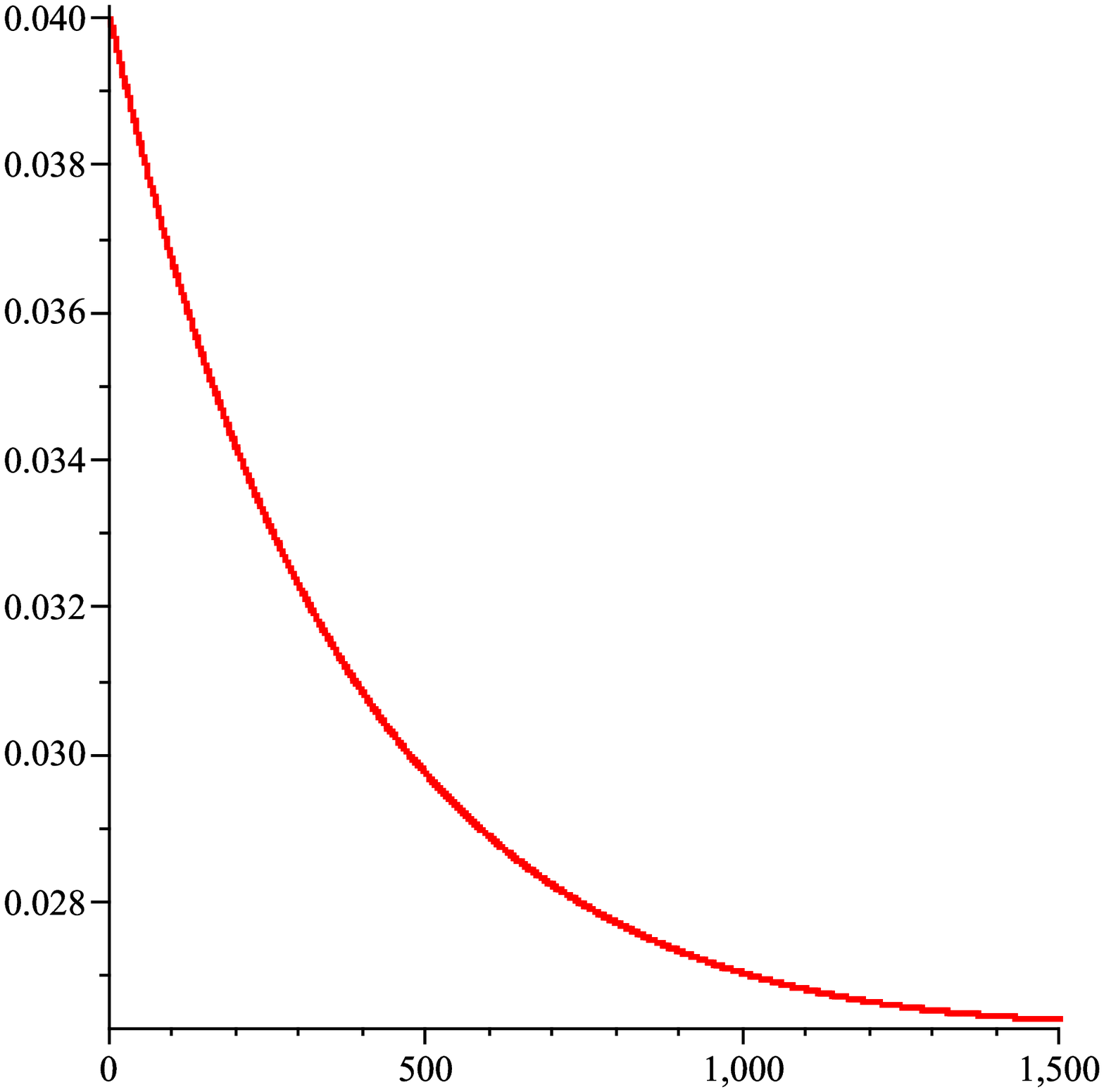} & \includegraphics[width=5cm]{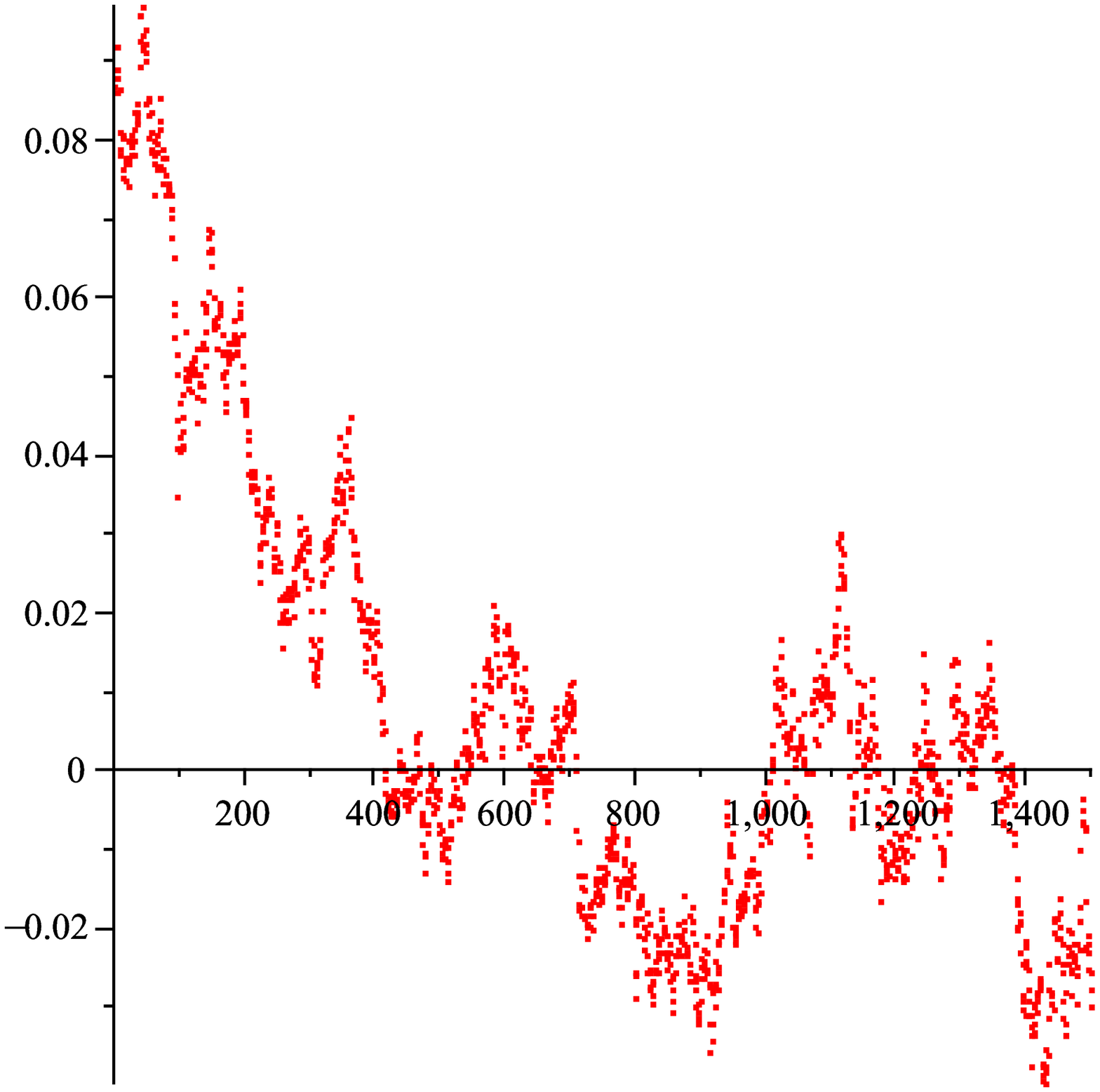}
\end{tabular}
\end{center}

Figures 5 and 6 represent the orbits $(q(nh), p(nh))$ and
$(q(\omega, nh), p(\omega, nh))$:

\begin{center}\begin{tabular}{cc}
\\ Fig 5. $(q(nh), p(nh))$ & Fig 6. $(q(\omega,
nh), p(\omega, nh))$\\
\includegraphics[width=5cm]{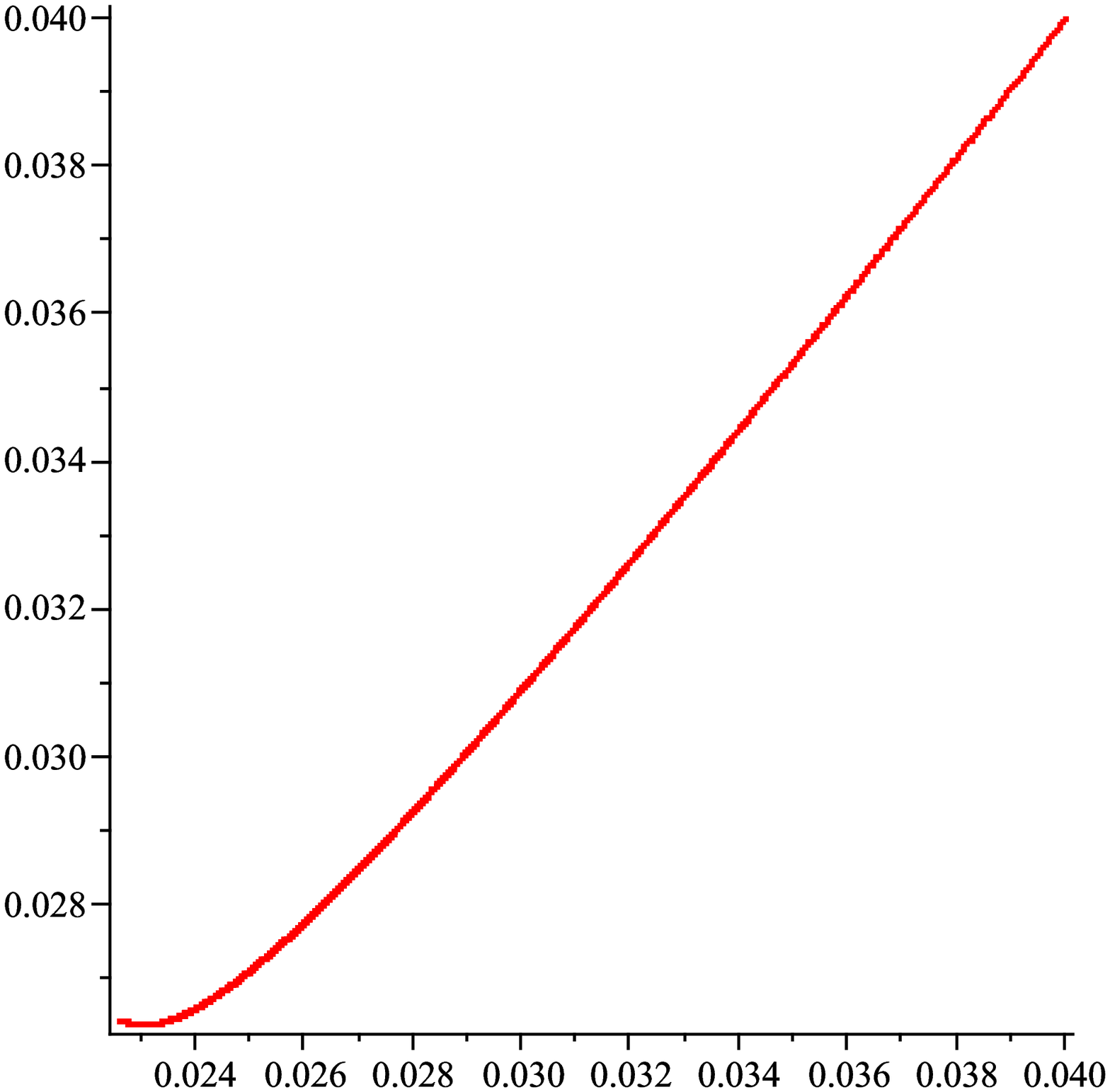} & \includegraphics[width=5cm]{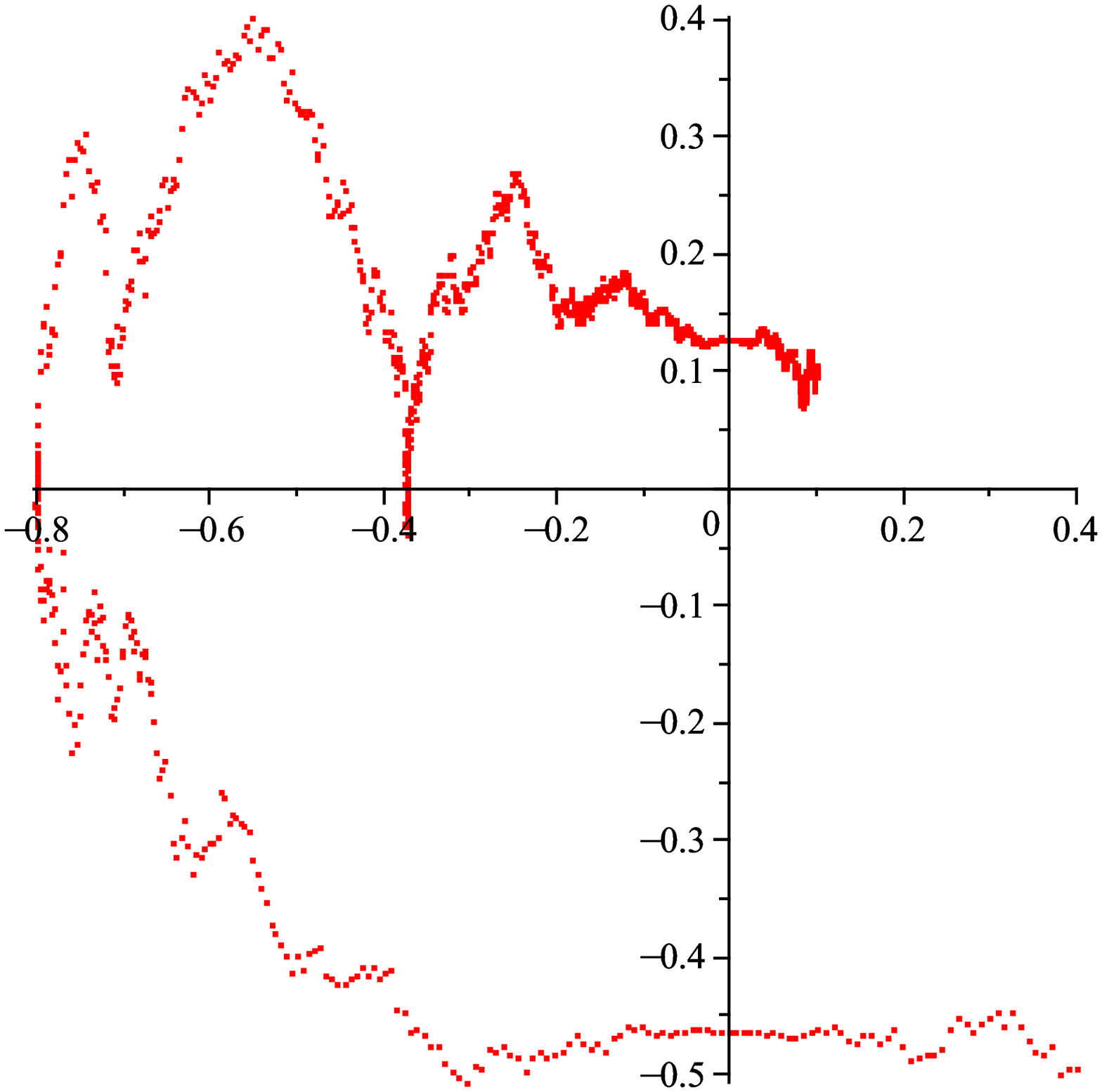}
\end{tabular}
\end{center}

\section{Stochastic generalized fractional HP}

\hspace{0.5cm} In this section a generalized variational principle
is introduced for a class of stochastic generalized fractional
Hamiltonian systems on  manifold. We use the generalized left
fractional Riemann-Liouville integral ([7]) defined as a mixture of
the fractal action from physics and the discounted action at rate
$\rho $. Let $f:\R\to\R$ be an integrable function, $\alpha
:\R\to\R$ a function of class $C^1$. The generalized left fractional
Riemman-Liouville integral is given by:

\begin{equation}
_{t_0}I^{\alpha}_tf(t)=\int_{t_0}^t\ds\frac{1}{\Gamma_1 (\alpha
(s-t))}f(s)(t-s)^{\alpha (s-t)-1}e^{\rho (s-t)}ds
\end{equation} where
\begin{equation*}
\Gamma_1 (\alpha (s-t))=\Gamma (\alpha (z))|z=s-t, \Gamma (\alpha
(z))=\int_0^{\infty}(s-t)^{\alpha (z)-1}e^{-(s-t)}ds,
\end{equation*} is the modified Euler Gamma function.

If $\alpha (z)=a=const.$, $0<a\leq 1$, $\rho =0$, from (15) we
obtain the left fractional Riemman-Liouville integral ([6], [7],
[8]). In fact, the generalized left fractional Riemann-Liouville
integral is a generalization of the single time Stieltyes integral
([13]).

In (15), $s$ is the intrinsec time and $t$ is the observer time,
$t\neq s$. Let $g:\R\to\R$ be the function:
\begin{equation}
g_t(s)=\ds\frac{1}{\Gamma_1 (\alpha (s-t))}e^{(\alpha
(s-t)-1)ln|t-s|+\rho (s-t)}, t\neq s.
\end{equation}

We consider ${\cal L}:J^1(\R, Q)\to\R$ and $\gamma :Q\to\R$. In the
stochastic context the HP principle states the following critical
point condition on $J^1(\R, Q)\oplus T^\ast (Q)$ for stochastic HP
generalized fractional action given by:

\begin{equation}
\begin{split}
& {\cal A}^{\alpha }(q,v,p, t)\!=\!\int_{t_0}^t[{\cal L}
(s,q(s),v(s))\!+\!<p(s),\ds\frac{dq}{ds}\!-\!v(s)>]g_t(s)ds\!+\!\\
& +\int_{t_0}^t\gamma (q(s))g_t(s)dw(s).
\end{split}
\end{equation}

The first integral in (17) is a Lebeques integral with respect to
$s$ and the second one is an It\^{o} integral.

Using Theorem 1, we get:

{\bf Theorem 2.} {\it If ${\cal L}:J^1(\R, Q)\to\R$ and $\gamma
:Q\to\R$ satisfy the hypothesis from Theorem 1, then almost
certainly a curve $c=(q,v,p)\in C(PQ)$ satisfies the stochastic HP
equations with intrinsec and observer times:

\begin{equation}
\begin{split}
& dq^i=v^{i}ds, \\
& dp_i=[\ds\frac{\pa {\cal L}}{\pa q^i}-p_i(\ds\frac{d\alpha
(s-t)}{ds}ln|t-s|+\ds\frac{\alpha (s-t)-1}{s-t})+\rho -\\
& -\ds\frac{1}{\Gamma_1 (\alpha (s-t))}\ds\frac{d\Gamma_1 (\alpha
(s-t))}{ds}]ds+\ds\frac{\partial \gamma (q)}{\partial q^i}dw(s),\\
& p_i=\ds\frac{\pa {\cal L}}{\pa v^i}, i=1..n, t\neq s.
\end{split}
\end{equation}}

From (18) we obtain:

(i) If $\alpha (z)=1$ and $\rho =0$ then equations (2) are obtained;

(ii) If $\alpha (z)=a=const.$, $0<a\leq 1$, $\rho =0$ then the
following relations are deduced from (18):

\begin{equation}
\begin{split}
& dq^i=v^{i}ds, \\
& dp_i=(\ds\frac{\pa {\cal L}}{\pa q^i}-p_i\ds\frac{a-1}{s-t})ds
+\ds\frac{\partial \gamma (q)}{\partial q^i}dw(s),\\
& p_i=\ds\frac{\pa {\cal L}}{\pa v^i}, i=1..n, t\neq s.
\end{split}
\end{equation}

The equations (19) represent the stochastic fractional equations.

If ${\cal L}:J^1(\R, Q)\to\R$ is hyperregular, using (18) the
following propositions hold:

{\bf Proposition 6}. (Stochastic generalized fractional Hamilton
equations.) {\it Equations (18) are equivalent with the equations:

\begin{equation}
\begin{split}
& dq^i=\ds\frac{\pa H}{\pa p_i}ds, \\
& dp_i=(-\ds\frac{\pa H}{\pa q^i}-p_ih(s,t))ds
+\ds\frac{\partial \gamma (q)}{\partial q^i}dw(s),\\
\end{split}
\end{equation} where
\begin{equation}
\begin{split}
& H=p_iv^i-{\cal L}(s,q,v), \\
& h(s,t)=\ds\frac{d\alpha (s-t)}{ds}ln|t-s|+\ds\frac{\alpha
(s-t)-1}{s-t}+\rho -\ds\frac{1}{\Gamma_1 (\alpha
(s-t))}\ds\frac{d\Gamma_1 (\alpha (s-t))}{ds}.
\end{split}
\end{equation}}
The equations (20) represent the generalized fractional Langevin
equations.

{\bf Proposition 7}. {\it If ${\cal L}=\ds\frac{1}{2}g_{ij}v^iv^j$,
where $g_{ij}$ are the components of a metric on the manifold $Q$,
then equation (18) takes the form:

\begin{equation}
\begin{split}
& dq^i=v^ids, \\
& dv^i=-(\Gamma ^i_{jk}v^jv^k-h(s,t)v^i)ds+g^{ij}\ds\frac{\pa \gamma
(q)}{\pa q^j}dw(s), i,j=1..n,
\end{split}
\end{equation} where $\Gamma^i_{jk}$ are Cristoffel coefficients associated to the considered metric
and $h(s,t)$ is given by (21). The equations (20) become:
\begin{equation}
\begin{split}
& dq^i=g^{ij}p_jds, \\
& dp_i=(\ds\frac{1}{2}\ds\frac{\pa g_{kl}}{\pa
q^i}p^lp^k-h(s,t)p_i)ds+\ds\frac{\pa \gamma (q)}{\pa q^i}dw(s),
i=1..n,
\end{split}
\end{equation} where $h(s,t)$ is given by (21).}

The equations (23) can be used for generalized fractional motion of
relativistic particle with noise ([7]).

{\bf Proposition 8}. {\it If ${\cal L}:J^1(\R, \R)\to\R$  is given
by:
\begin{equation}
{\cal L}(q,v)=\ds\frac{1}{2}v^2-V(q)
\end{equation} where $V:\R\to\R$ and $\gamma :\R\to\R$, then
equations (23) are given by}:

\begin{equation}
\begin{split}
& dq=pds, \\
& dp=(-\ds\frac{dV}{dq}-h(s,t)p)ds+\ds\frac{d \gamma (q)}{d q}dw(s).
\end{split}
\end{equation}

If $V(q)=cos(q)$, $\gamma (q)=sin(q)$, the Euler scheme for (25) is
given by:
\begin{equation*}
\begin{split}
& q(n+1)=q(n)+kp(n)) \\
& p(n+1)=p(n)+k(sin(q(n))-h(nk, t))+cos(q(n))G(n), n=0..N-1,
\end{split}
\end{equation*} where $T>0$, $k=\ds\frac{T}{N}$, $N>0$,
$G(n)=w((n+1)h)-w(nh)$ and

\begin{equation*}
\begin{split}
& h(nk, t)=\alpha ^1(nk, t)ln|t-nk|+\ds\frac{\alpha (nk-t)-1}{nk-t}+\rho-\\
&-\ds\frac{1}{\Gamma_1 (\alpha (nk-t))}\Gamma ^1(\alpha (nk-t)), \\
& \alpha^1(s-t)=\ds\frac{d\alpha (s-t)}{ds}, \Gamma^1(\alpha
(s-t))=\ds\frac{d\Gamma_1 (\alpha (s-t))}{ds}.
\end{split}
\end{equation*}

For $\alpha (s-t)=a$, $a=0.6$, $k=0.001$, $t=0.8$ using Maple 13 the
orbit $(n,p(nk))$ is represented in Figure 7 and the orbit
$(n,p(\omega , nk))$ is represented in Figure 8:

\begin{center}\begin{tabular}{cc}
\\ Fig 7. $(n,p(nk))$ & Fig 8. $(n,p(\omega , nk))$\\
\includegraphics[width=5cm]{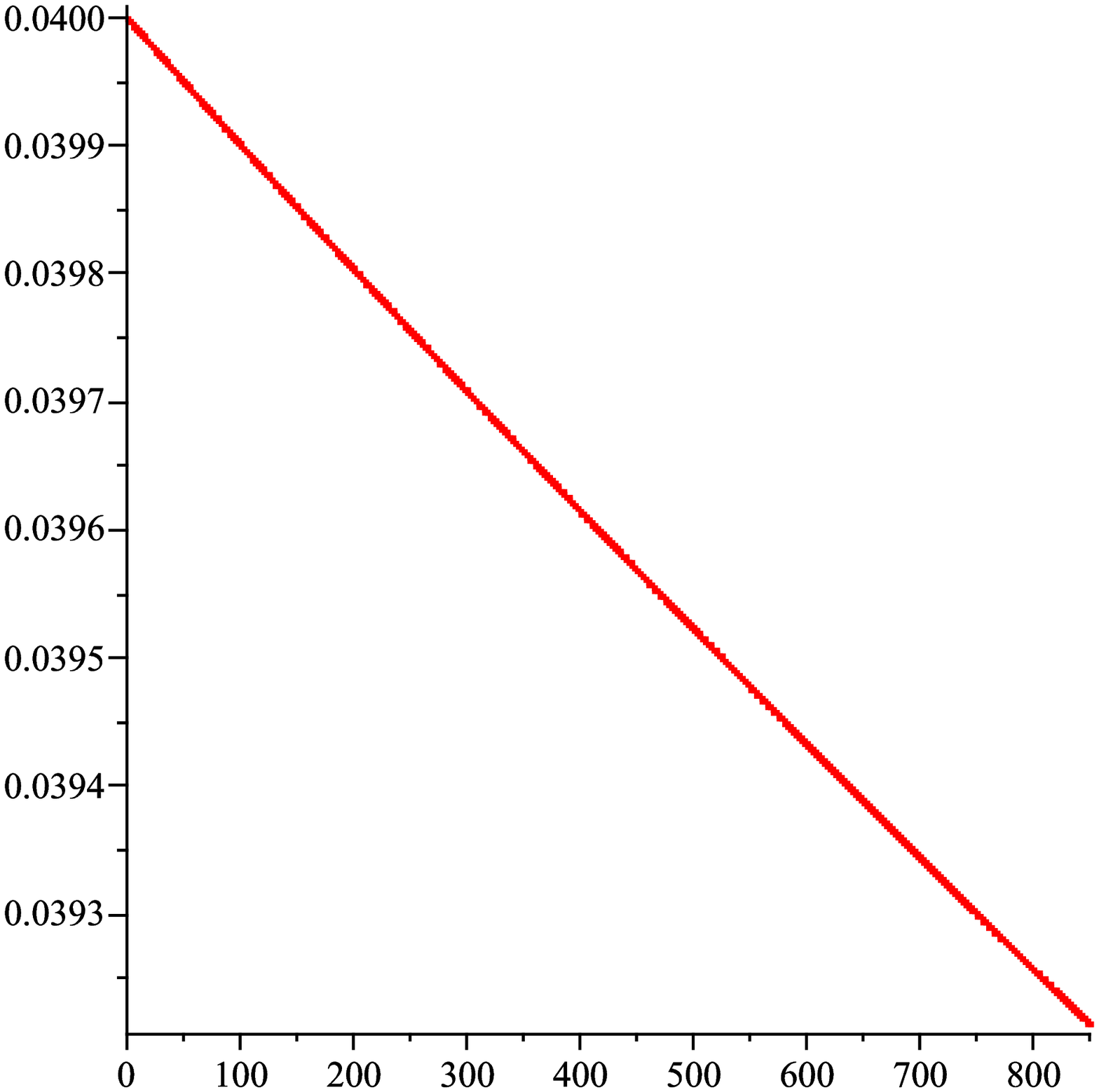} & \includegraphics[width=5cm]{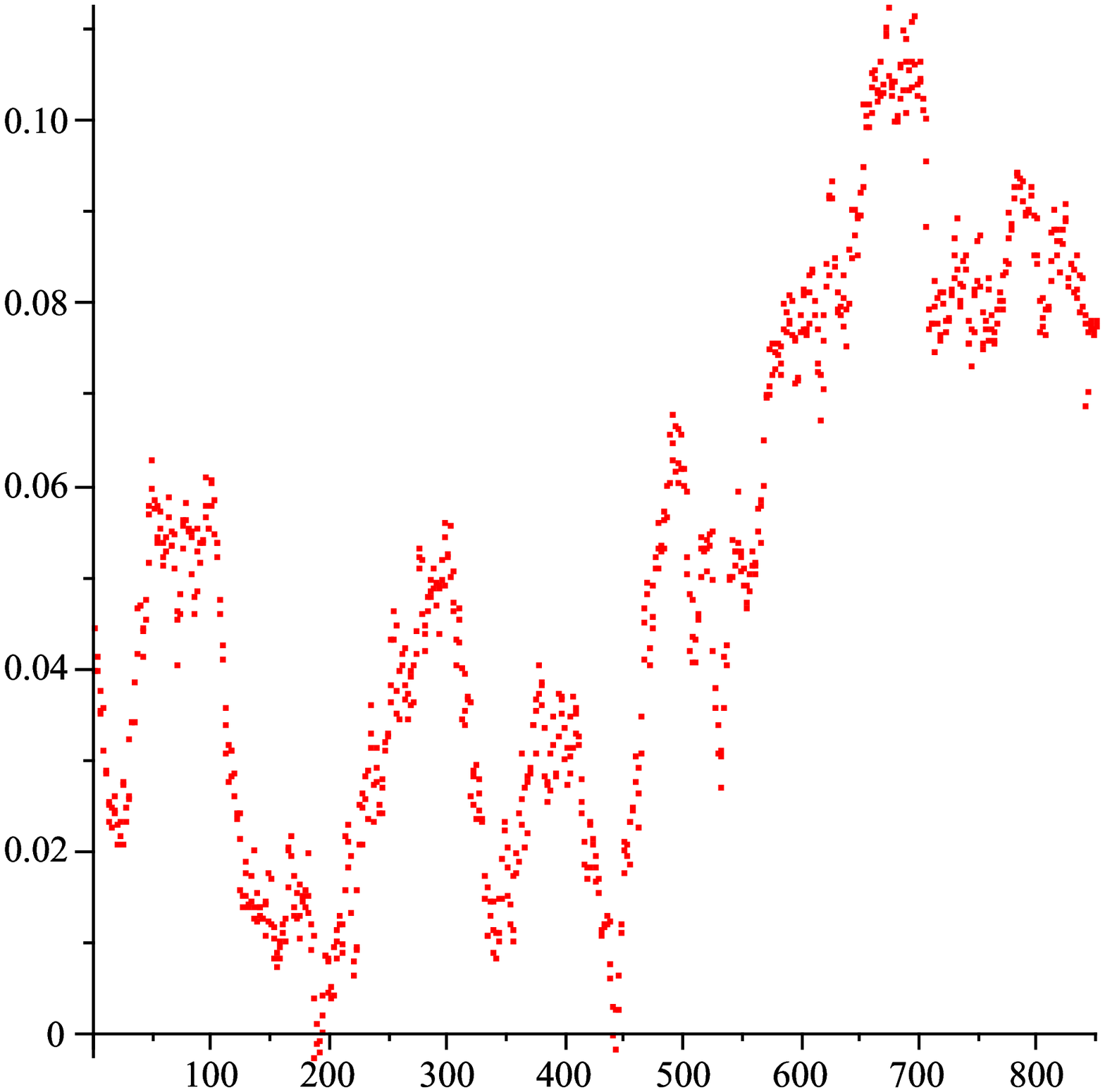}
\end{tabular}
\end{center}

Figure 9 displays the orbit $(q(nk), p(nk))$ and Figure 10 displays
the orbit $(q(\omega, nk), p(\omega, nk))$

\begin{center}\begin{tabular}{cc}
\\ Fig 9. $(q(nk), p(nk))$ & Fig 10. $(q(\omega, nk), p(\omega, nk))$\\
\includegraphics[width=5cm]{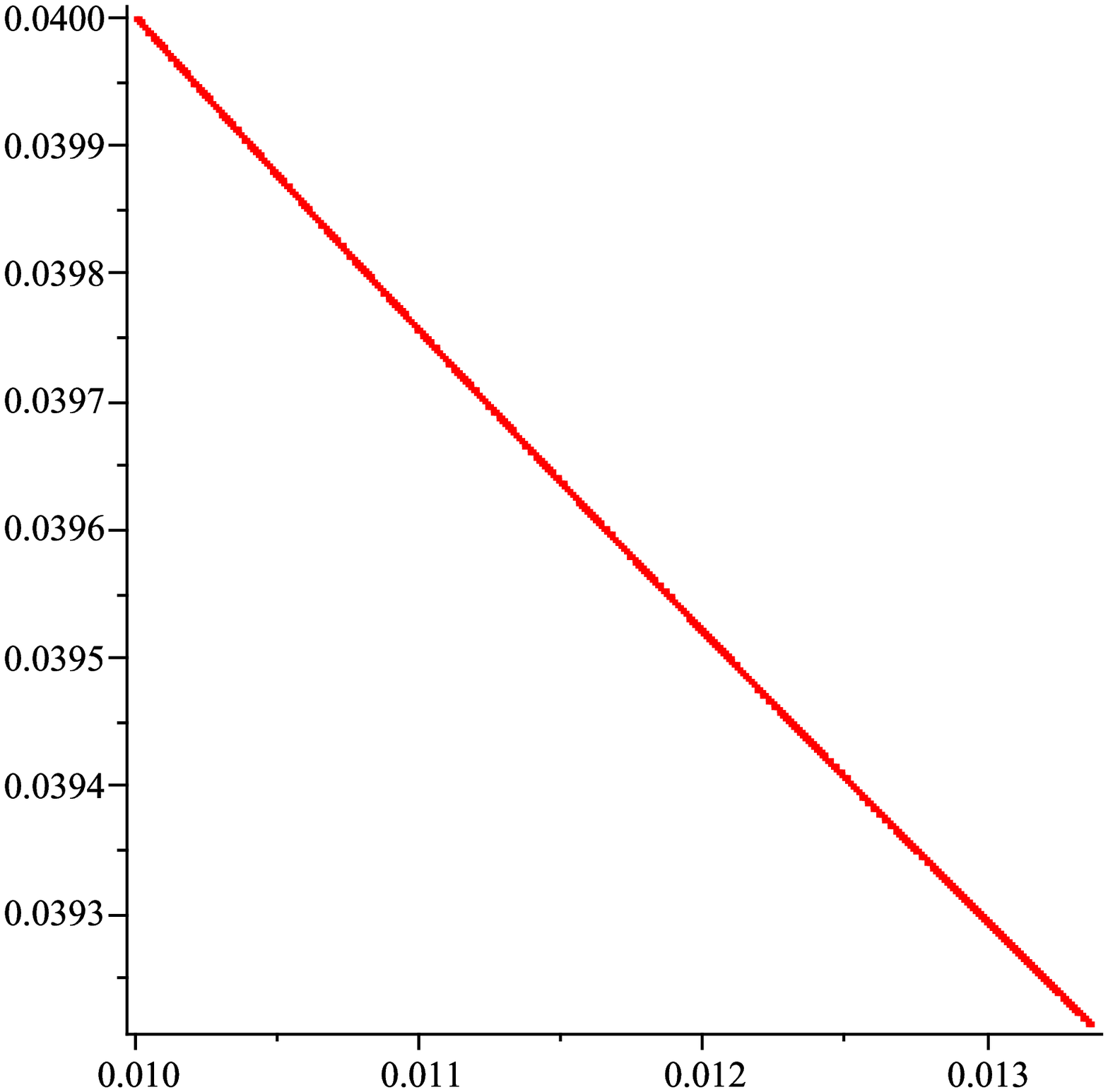} & \includegraphics[width=5cm]{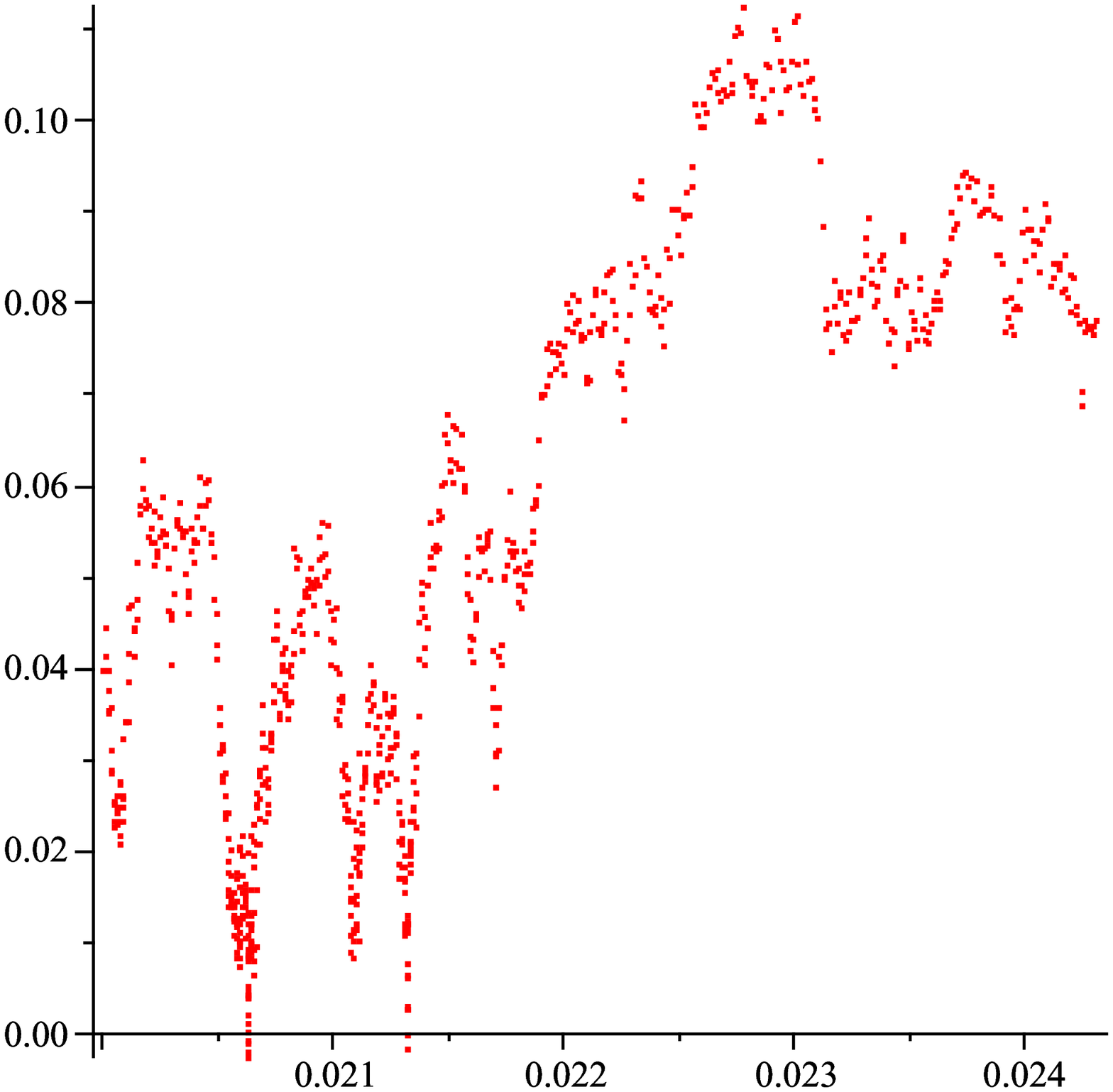}
\end{tabular}
\end{center}

\section{Conclusions}

\hspace{0.5cm} In this paper we have described the stochastic
generalized fractional HP principle, using the classical stochastic
HP principle [3]. Using a hyperregular Lagrange function,
Langevin-type generalized fractional equations were illustrated. We
have done the numerical simulations for the case of a Lagrange
function defined on $\R^2$. In our future papers we will study the
stochastic stability of the obtained equations al well as the
description of the credibility generalized fractional HP principle.

{}

\end{document}